\documentclass[12pt]{article}
\usepackage{amsmath}

\usepackage{amsthm} 
\usepackage{amscd}
\theoremstyle{plain} 
\newtheorem{theorem}{Theorem}[section]
\newtheorem{corollary}[theorem]{Corollary}

\newtheorem{lemma}[theorem]{Lemma}

\newtheorem{definition}[theorem]{Definition}
\numberwithin{equation}{section}

\newcommand{\al}{\alpha}
\newcommand{\be}{\beta}
\newcommand{\br}{\mathbf {R}}
\newcommand{\om}{\omega}
\newcommand{\ot}{\otimes}
\newcommand{\ld}{\ldots}

\newcommand{\g} {\mathsf g }

\begin{document}
\title{Rigidity of Secondary Characteristic Classes}
\author{Jerry M. Lodder}
\date{}
\maketitle

\section{Introduction}

In this paper we study the variability and rigidity of secondary
characteristic classes which arise from flat connections $\theta$ on a
differentiable manifold $M$.  In particular we consider $\theta$ as a
Lie-algebra valued one-form on $M$, and study the characteristic map
$$  \phi_{\theta}: H^*_{\text{Lie}}(\g ) \to H^*_{\text{dR}}(M), $$
where $H^*_{\text{Lie}}$ denotes Lie algebra cohomology (for a Lie
algebra $\g$), and $H^*_{\text{dR}}$ denotes de Rham cohomology.  An
element $\al \in H^*_{\text{Lie}}(\g )$ is called variable if there exists
a one-parameter family of flat connections $\theta_t$ with
$$  \frac{d}{dt} \Big[ \phi_{\theta_t} (\al) \, 
	\Big{\vert}_{t=0} \neq 0,  $$
otherwise $\al$ is called rigid.  For example, the universal
Godbillon-Vey invariant for codimension $k$ foliations is known to be
a variable class.  Letting $HL^*$ denote Loday's Leibniz
cohomology \cite{Ens.Math.}
\cite{overview} \cite{LP}, we prove that if $HL^n (\g )= 0$
for $n \geq 1$, then all classes in $H^*_{\text{Lie}}(\g )$ are rigid.
This result imparts a geometric meaning to Leibniz cohomology.

Moreover, in the case of codimension one foliations (with trivial
normal bundle), there is a flat $W_1$-valued connection on $M$, where
$W_1$ is the Lie algebra of formal vector fields on $\br^1$.  
From \cite{Lodder},
$$  HL^*(W_1 ) \simeq \Lambda(\al) \ot T(\zeta), $$
where $\Lambda (\al)$ is the exterior algebra on the Godbillon-Vey
invariant (in dimension three) and $T(\zeta)$ denotes the tensor
algebra on a four-dimensional class.  When $M$ is provided with a
one-parameter family of such foliations, we compute the image of a
characteristic map
$$  HL^4 (W_1 ) \to H^4 _{\text{dR}}(M).  $$
The image in de Rham cohomology is independent of the choices made
when constructing the $W_1$ connections, and involves the time
derivative (derivative with respect to $t$) of the Bott connection.

\section{The Characteristic Map and Rigidity}

Let $M$ be a differentiable ($C^{\infty}$) manifold with flat
connection $\theta$.  We consider the formulation of $\theta$ as a
Lie-algebra valued one-form
$$  \theta : TM \to \g , $$
where $TM$ denotes the tangent bundle of $M$ and $\g$ is a Lie
algebra.  This section describes a characteristic map
$$ \phi_{\theta} : H^*_{\text{Lie}} (\g ; \, \br) \to
	H^*_{\text{dR}}(M) $$
from Lie algebra cohomology (with $\br$ coefficients) to the de Rham
cohomology of $M$.  In the case of a topological Lie algebra, then 
$H^*_{\text{Lie}}$ is understood as continuous cohomology, computed
using continuous cochains.  

Let $\Omega^k (\g ; \, \br )$ be the $\br$-vector space of
skew-symmetric (continuous) cochains
$$  \al : \g^{\ot k} \to \br .  $$
For comparison with Leibniz cohomology (and to establish our sign
conventions), we write the coboundary for
Lie algebra cohomology
$$  d: \Omega^k (\g; \, \br) \to \Omega^{k+1}(\g; \, \br) $$
as 
\begin{equation} \label{2.1}
\begin{split}
& d(\al)(g_1 \ot g_2 \ot \, \ld \, \ot g_{k+1}) =  \\
& \sum_{1 \leq i < j \leq k} (-1)^{j+1} \al (g_1 \ot \, \ld \, g_{i-1}
\ot [g_i, \, g_j]  \ot g_{i+1} \, \ld \, \hat{g}_j \, \ld \ot g_{k+1}),
\end{split}
\end{equation}
where each $g_i \in \g$.  For a differentiable manifold $M$, let
$$  \Omega^k (M) := \Omega^k(M; \, \br) $$
denote the $\br$-vector space of $k$-forms on $M$.  Then the de Rham
coboundary 
\begin{align*}
&  d: \Omega^k (M) \to \Omega^{k+1}(M) \\
& \om \mapsto d \om
\end{align*}
has a global formulation as 
\begin{equation} \label{2.2}
\begin{split}
  & d \om  (X_1 \ot  X_2 \ot \, \ldots \, \ot X_{k+1}) =  \\
  & \sum_{i=1}^{k+1} (-1)^{i+1} X_i \big( \om (X_1 \ot \, \ldots \,
   \hat{X}_i \,  \ldots \ot X_{k+1}) \big) +  \\
  & \sum_{1 \leq i < j \leq k+1} (-1)^{j+1} \om \big( X_1 \ot \, \ldots
   \, \ot X_{i-1}\ot   [X_i, \, X_j]\ot X_{i+1}\ot \, \\  
  & \hskip2.5in \ldots \, 
	\hat{X}_j \, \ldots \, \ot X_{k+1} \big), 
\end{split}
\end{equation}
where each $X_i \in \chi(M)$, the Lie algebra of smooth ($C^{\infty}$)
vector fields on $M$.  

With the above sign conventions for the coboundary maps, the 
connection $\theta \in \Omega^1(M; \, \g )$ is flat if and only if the
Maurer-Cartan equation holds,
\begin{equation} \label{2.3}
  d \theta = - \frac{1}{2} \, [\theta, \, \theta], 
\end{equation}
where $ d \theta \in \Omega^2 (M; \, \g)$ is given by
$$  (d \theta)(X_1 \ot X_2) = X_1 (\theta(X_2)) - X_2(\theta (X_1))
	- \theta([X_1, \, X_2]).  $$
Recall that for a differentiable function $f : M \to \g$, $X(f)$
denotes the partial derivatives of the component functions of $f$
with respect to $X$.  If $e_1$, $e_2$, $e_3$, $\ld$ is a basis for
$\g$, then
$$ f = f_1 e_1 + f_2 e_2 + f_3 e_3 + \, \cdots , $$
where each $f_i : M \to \br$ is differentiable.  Thus,
$$  X(f) = X(f_1)e_1 + X(f_2)e_2 + X(f_3)e_3 + \, \cdots .$$
Also, the symbol $[\theta, \, \theta]$ denotes the composition 
$$  \Omega^1(M;\, \g) \ot \Omega^1(M;\, \g) 
	\overset {\theta \wedge \theta}{\longrightarrow}
	\Omega^2(M; \, \g \ot \g) 
	\overset{[ \; , \; ]}{\longrightarrow} \Omega^2(M; \, \g)$$
given by
\begin{align*}
[\theta , \, \theta] (X_1 \ot X_2)\, & = [ \ \, , \;] \circ \big(
  \theta(X_1) \ot \theta(X_2) - \theta(X_2) \ot \theta(X_1)\big) \\
& = 2 [\theta (X_1), \, \theta(X_2)].
\end{align*}

Consider the following map \cite[p.\ 234]{Fuks}
$$  \phi: \Omega^*(\g ; \, \br) \to \Omega^*(M; \, \br)$$
given on a cochain $\al \in \Omega^k(\g; \, \br)$ by
$$
 \phi(\al)(X_1 \ot X_2 \ot \, \ld \, \ot X_k) = 
 \al(\theta(X_1) \ot \theta(X_2) \ot \, \ld \, \ot \theta(X_k)).
$$
For the case $k=1$, it can be easily seen from (2.3) that $\phi$ is a
map of cochain complexes.  In particular,
\begin{align*}
d(\phi(\al & )) (X_1 \ot X_2)  \\
& = X_1 (\al (\theta (X_2))) - X_2(\al (\theta (X_1))) - 
	\al(\theta( [X_1, \, X_2])) \\
& = \al \Big( X_1 (\theta(X_2)) - X_2(\theta(X_1)) - 
	\theta([X_1, \, X_2])\Big) \\
& = \al(- [\theta(X_1), \, \theta(X_2)]) \\
& = \phi(d\al)(X_1 \ot X_2).
\end{align*}
In \cite[p.\ 235]{Fuks} it is proven that in general 
$$  d \phi = \phi d . $$
The induced map
$$ \phi_{\theta} : H^*_{\text{Lie}}(\g ) \to H^*_{\text{dR}}(M) $$
is called the characteristic map.

A specific example of a flat connection studied in this paper arises
from the theory of foliations.  Let $\cal F$ be a $C^{\infty}$
codimension one foliation on $M$ with trivial normal bundle.  Given a
choice of a trivialization, then a determining one-form $\om_0$ is
defined for the foliation by $\om_0(v) = 0$ is $v$ is tangent to a
leaf and $\om_0(\eta) = 1$ if $\eta$ is a unit vector of the normal
bundle with positive orientation.  Letting $d$ denote the de Rham
coboundary, a sequence of one-forms $\om_1$, $\om_2$, $\om_3$, $\ld$
can be defined so that \cite{Godbillon}
\begin{equation} \label{2.4}
\begin{split}
& d \om_0 = \om_0 \wedge \om_1  \\
& d \om_1 = \om_0 \wedge \om_2  \\
& d \om_2 = \om_0 \wedge \om_3 + \om_1 \wedge \om_2 \\
& d \om_k = \sum_{i=0}^{\textstyle{[\frac{k}{2}}]} \, \, \frac{k-2i+1}{k+1}
 \, \binom{k+1}{i} \, (\om_i \wedge \om_{k+1-i}).
\end{split}
\end{equation}
Consider the topological Lie algebra of formal vector fields
$$ W_1 = \Big\{ \, \sum_{k=0}^{\infty} c_k\, x^k \, \frac{d}{dx} 
	\ \ \Big{\vert}\ \ c_k \in \br \, \Big\}  $$
in the $\cal M$-adic topology, where $\cal M$ is the maximal ideal of
$R[[x]]$ given by those series with zero constant term.  Then a
$W_1$-valued one-form is defined on $M$ by 
\begin{equation} \label{2.5}
 \theta_{\cal F}(v) = \sum_{k=0}^{\infty} \om_k(v) \, 
	\frac{x^k}{k!} \, \frac{d}{dx},
\end{equation}
where $v \in TM$.  From (2.4), it can be proven \cite[p.\ 231]{Fuks}
that $\theta_{\cal F}$ is a flat connection on $M$.  The resulting
homomorphism 
\begin{equation} \label{1.6}
\phi_{\theta_{\cal F}}: H^*_{\text{Lie}}(W_1) \to H^*_{\text{dR}}(M)
\end{equation} 
is the classical characteristic map in foliation theory.  (A similar
construction exists for foliations of any codimension with trivial
normal bundle.)  

We wish to study a one-parameter variation of a flat structure 
\begin{equation} \label{2.7}
 	\theta_t : TM \to \g, \ \ \ t \in \br, \ \ \ \theta_0 = \theta,
\end{equation}
which depends smoothly on the parameter $t$.  Such a structure may
arise from a one-parameter variation of a foliation.
\begin{definition} \label{2.8}
A Lie algebra cohomology class $\al \in H^*_{\text{Lie}}(\g)$ is called
variable (for $\theta$) if there exists a family $\theta_t$ such that 
$$  \frac{d}{dt} \Big[ \phi_{\theta_t} (\al) \, 
	\Big{\vert}_{t=0} \neq 0.  $$
Otherwise, $\al$ is called rigid.
\end{definition}
By work of Thurston, the universal Godbillon-Vey invariant $\al \in 
H^3_{\text{Lie}}(W_1)$ is a variable class \cite{Thurston}.  One of
the goals of 
this paper is to prove that if the Leibniz cohomology of $\g$
vanishes, i.e., $HL^n(\g ) = 0$, $n \geq 1$, then all characteristic
classes in $H^*_{\text{Lie}}(\g )$ are rigid.  In the remainder of
this section we restate the definition of rigidity in terms of a known
condition concerning $H^*_{\text{Lie}}(\g ;\, \g')$, the
Lie algebra cohomology of $\g$ with coefficients in the
coadjoint representation 
$$ \g' = {\text{Hom}}^c_{\br}(\g, \, \br),  $$
($c$ denotes continuous maps).  

First introduce the current algebra $\tilde{\g} = C^{\infty}(\br, \,
\g )$ of differentiable maps from $\br$ to $\g$.  Then given
$\theta_t$ as in (2.7), there is a flat $\tilde{\g}$ connection on $M$
\begin{equation} \label{2.9}
\begin{split}
	& \Theta : TM \to \tilde{\g}  \\
	& \Theta(v)(t) = \theta_t(v), \ \ \ v \in TM,
\end{split}
\end{equation}
and a characteristic map
$$  \phi_{\Theta} : H^*_{\text{Lie}}( \tilde{\g}) \to H^*_{\text{dR}}(M).$$
Using an idea of D. Fuks \cite{Fuks}, 
define a ``time derivative'' map on cochains
$$  D: \Omega^q (\g) \to \Omega^q(\tilde{\g})  $$
by 
\begin{align*}
& D(\al)(\varphi_1, \, \varphi_2, \, \ld , \, \varphi_q) =  \\
& \sum_{i=1}^{q} \al \Big(\varphi_1(0), \, \ld , \,
 \varphi_{i-1}(0), \,  \frac{d}{dt}[\varphi_i(t)\,\big{\vert}_{t=0}\, , \,
 \varphi_{i+1}(0), \, \ld , \, \varphi_q (0) \Big),
\end{align*}
where $\al \in \Omega^q (\g )$, $(\varphi_1, \, \varphi_2, \, \ld, \,
\varphi_q ) \in (\tilde{\g})^{\ot q}$.  Then $D$ is a map of cochain
complexes, and there is an induced map
\begin{equation} \label{1.10}
  D^* : H^*_{\text{Lie}}(\g ) \to H^*_{\text{Lie}}(\tilde{\g}).  
\end{equation}
It follows that given $\al \in H^*_{\text{Lie}}(\g )$, we have
\begin{equation} \label{2.11}
 \phi_{\Theta} \circ D^* (\al) = \frac{d}{dt} \big[ \phi_{\theta_t}
	(\al ) \big{\vert}_{t=0}\, .
\end{equation}

Recall that $\g' = \text{Hom}^c_{\br}(\g , \, \br)$ is a left
$\g$-module with 
$$  (g \gamma)(h) = \gamma([h, \, g]), $$
where $g, \, h \in \g$ and $\gamma \in \g'$.  Then $D^*$ can be
factored as $\varPhi^* \circ V^*$ \cite[p.\ 244]{Fuks}, where
\begin{equation} \label{2.12}
\begin{split}
& V^* : H^q_{\text{Lie}}( \g) \to H^{q-1}(\g; \, \g'), \ \ q \geq 1, \\
& \varPhi^* : H^{q-1}(\g; \, \g') \to H^q_{\text{Lie}}(\tilde{\g}),
  \ \ q \geq 1, 
\end{split}
\end{equation}
are induced by
\begin{equation} \label{2.13}
\begin{split}
& \text{Var}: \Omega^q(\g; \, \br) \to \Omega^{q-1}(\g; \, \g')    \\
& \varPhi : \Omega^{q-1}(\g; \, \g') \to \Omega^q(\tilde{\g}; \, \br) \\
& (\text{Var})(\al)(g_1, \, g_2, \, \ld, \, g_{q-1})(g_0) = 
  (-1)^{q-1} \al(g_0, \, g_1, \, \ld, \, g_{q-1})  \\
& \varPhi (\gamma)(\varphi_1, \, \varphi_2, \, \ld, \, \varphi_q) = \\
& \ \ \ \ \ 
  \sum_{i=1}^q (-1)^{q-i} \, \gamma \big( \varphi_1(0), \, \ld, \,
  \hat{\varphi}_i(0), \, \ld, \, \varphi_q(0) \big) \big( \varphi'_i(0)\big),
\end{split}
\end{equation}
where $\al \in \Omega^q(\g; \, \br)$, $g_i \in \g$,
$\gamma \in \Omega^{q-1}(\g; \, \g')$, $\varphi_i \in \tilde{\g}$.
We then have a commutative diagram
\begin{equation} \label{2.14}
\CD  
	H^*_{\text{Lie}}(\g) @>D^*>> H^*_{\text{Lie}}(\tilde{\g}) \\
	@V{V^*}VV   @AA{\varPhi^*}A  \\
	H^{*-1}_{\text{Lie}}(\g; \, \g') @>=>>
		H^{*-1}_{\text{Lie}}(\g; \, \g') 
\endCD
\end{equation}
\begin{lemma} \label{2.15}
If $H^{n-1}_{\text{Lie}}(\g; \, \g') = 0$ for $n \geq 1$, then all
characteristic classes in $H^*_{\text{Lie}}(\g )$ are rigid.
\end{lemma}
\begin{proof}
This follows from equation \eqref{2.11}, diagram \eqref{2.14} and the 
definition of rigidity (definition \eqref{2.8}).
\end{proof}

In the next section we prove that if $HL^n(\g) = 0$ for $n \geq 1$,
then $H^{n-1}_{\text{Lie}}(\g; \, \g') = 0$ for $n \geq 1$.

\section{Leibniz Cohomology}

Still considering $\g$ to be a Lie algebra (over $\br$), recall that 
the Leibniz cohomology of $\g$ with trivial coefficients,
$$  HL^*(\g; \, \br) := HL^*(\g ),  $$
is the homology of the cochain complex \cite{LP}
\begin{equation} \label{3.1}
  \br \overset{0}{\rightarrow} C^1(\g) \overset{d}{\rightarrow}
  C^2(\g) \rightarrow \, \cdots \, \rightarrow C^k(\g) \overset{d}{\rightarrow}
  C^{k+1}(\g) \rightarrow \,  \cdots \, ,
\end{equation}
where $C^k(\g ) = {\text{Hom}}^c_{\br}( \g^{\ot k}, \, \br)$, and for
$\al \in C^k(\g)$, $d\al$ is given in equation \eqref{2.1}.  Keep in mind that
for Leibniz cohomology, the cochains are not necessarily
skew-symmetric.  

In this section we prove the following:
\begin{theorem} \label{3.2}
If $HL^n(\g; \, \br) = 0$ for $n \geq 1$, then
$H^{n-1}_{\text{Lie}}(\g; \, \g') = 0$ for $n \geq 1$, where $\g' =
{\text{Hom}}^c_{\br}(\g; \, \br)$.  
\end{theorem}
The proof involves a spectral sequence similar to the Pirashvili
spectral sequence \cite{Pirashvili}, except tailored to the specific
algebraic relation between $HL^*(\g)$ and $H^{*-1}_{\text{Lie}}(\g; \,
\g')$.  Recall that the projection to the exterior power
$$  \g^{\ot q} \to \g^{\wedge q}  $$
induces a homomorphism
$$  H^*_{\text{Lie}}(\g) \to HL^*(\g).  $$
Letting $C^*_{\text{rel}}(\g)[2] = C^*(\g)/ \Omega^*(\g)$,
we have a long exact sequence
$$ \cdots \, \rightarrow H^q_{\text{Lie}}(\g ) \rightarrow HL^q(\g)
   \rightarrow H^{q-2}_{\text{rel}}(\g) \rightarrow
   H^{q+1}_{\text{Lie}}(\g ) \rightarrow \, \cdots \, . $$
The Pirashvili spectral sequence arises from a filtration of
$C^*_{\text{rel}}(\g)[2]$ and converges to $H^*_{\text{rel}}(\g)$.

Consider now the map of cochain complexes
$$  i: \Omega^{q-1}(\g; \, \g') \to C^q({\g})  $$
given by
$$  \big( i(\be )\big) (g_0 \ot g_1 \ot \, \ld \, \ot g_{q-1})
   = (-1)^{q-1} \, \be(g_1 \ot g_2 \ot \, \ld \, \ot g_{q-1})(g_0), $$
where $\be \in \Omega^{q-1}(\g; \, \g')$ and $g_i \in \g$ for 
$i = 0, \, 1, \, 2, \, \ld, \, q-1$.  Letting
$$  C^*_{RG}(\g)[2] = C^*(\g)/i[\Omega^{*-1}(\g; \, \g')],  $$
we also have a long exact sequence
\begin{equation} \label{3.3}
\begin{split}
& 0 \rightarrow H^0_{\text{Lie}}(\g; \, \g') \rightarrow HL^1(\g)
  \rightarrow 0 \rightarrow  \\
& H^1_{\text{Lie}}(\g; \, \g') \rightarrow
  HL^2(\g) \rightarrow H^0_{RG}(\g) \rightarrow  
  H^2_{\text{Lie}}(\g; \,  \g')\rightarrow  \, \\
&  \cdots \, \rightarrow H^{q-1}_{\text{Lie}}(\g; \, \g') \rightarrow
  HL^q(\g) \rightarrow H^{q-2}_{RG}(\g) \rightarrow H^q_{\text{Lie}}(\g; \, \g')
  \rightarrow \, \cdots \, .
\end{split}
\end{equation}
The filtration for the Pirashvili spectral sequence \cite{Lodder}
\cite{Pirashvili} can be immediately applied to yield a decreasing
filtration $\{ F^s_* \}_{s \geq 0}$ for $C^*_{RG}(\g)[2]$.  We use the
same grading as in \cite{Lodder} \cite{Pirashvili}, which becomes
$F^0_* = C^*_{RG}[2]$, and for $s \geq 1$,
\begin{align*}
&  F^s_* = A/B  \\
&  A = \{ \, f \in C^*(\g) \ | \ f \ {\text{is skew-symmetric in the
last $(s+1)$ tensor factors}} \, \}  \\
&  B = i [ \Omega^{*-1}(\g; \, \g')].
\end{align*}
Then as in \cite{Lodder}, each $F^s_*$ is a subcomplex of
$C^*_{RG}(\g)$, and
$$  F^0_* \supseteq F^1_* \supseteq F^2_* 
    \supseteq \, \ld \, \supseteq F^s_* \supseteq F^{s+1}_*
    \supseteq \, \ld \, .  $$

To identify the $E^2$ term of the resulting spectral sequence,
consider coker(Var), where Var is defined in equation \eqref{2.13}.
Letting 
$$  CR^*(\g)[1] = \Omega^{*-1}(\g; \, \g')/{\text{Var}}[\Omega^*(\g)], $$
there is a short exact sequence
\begin{equation} \label{3.4}
  0 \rightarrow \Omega^*(\g) \overset{\text{Var}}{\longrightarrow}
  \Omega^{*-1}(\g; \, \g') \rightarrow CR^*(\g)[1] \rightarrow 0,  
\end{equation}
and an associated long exact sequence
$$ \cdots \, \rightarrow H^q_{\text{Lie}}(\g) \rightarrow
  H^{q-1}_{\text{Lie}}(\g; \g') \rightarrow HR^{q-2}(\g) \rightarrow
  H^{q+1}_{\text{Lie}}(\g ) \rightarrow \, \cdots \, .  $$

\begin{theorem} \label{3.5}
The filtration $F^s_*$ of $C^*_{RG}(\g)[2]$ yields a spectral sequence
converging to $H^*_{RG}(\g)$ with
\begin{align*}
&  E^{s,0}_2 = 0, \ \ s = 0, \, 1, \, 2, \, \ld,  \\
&  E^{s,n}_2 \simeq HL^n(\g) \hat{\ot} HR^s(\g), \ \ n = 1, \, 2, \, 3,
  \, \ld, \ \ s = 0, \, 1, \, 2, \, \ld , 
\end{align*}
where $\hat{\ot}$ denotes the completed tensor product.
\end{theorem}
\begin{proof}
The proof follows from the identification of the $E^2$ term in
\cite{Lodder} or \cite{Pirashvili}.  Also, note that
$$  E^{s,0}_0 = (F^s/F^{s+1})_0, \ \ \ F^s_0 = A/B,  $$
where
\begin{align*}
&  A = \{ \, f \in C^{s+2}(\g) \ | \ f \ {\text{is alternating in the
	last $(s+1)$ factors}} \, \}  \\
&  B = i[\Omega^{s+1}(\g; \g')].  
\end{align*}
Then $A = B$, $F^s_0 = 0$, and $E^{s,0}_0 = 0$.  
\end{proof}

\begin{theorem} \label{3.6}
If $HL^n(\g) = 0$ for $n \geq 1$, then $H^{n-1}_{\text{Lie}}(\g; \,
\g') = 0$ for $n \geq 1$.
\end{theorem}
\begin{proof}
If $HL^n(\g) = 0$ for $n \geq 1$, then from theorem \eqref{3.5}, the
$E_2$ term for the spectral sequence converging to $H^*_{RG}(\g)$ is
zero.  Thus, $H^n_{RG}(\g) = 0$ for $n \geq 0$.  The result now
follows from long exact sequence \eqref{3.3}.
\end{proof}

\begin{theorem} \label{3.7}
If $HL^n(\g) = 0$ for $n \geq 1$, and $\theta_t$ is a one-parameter
family of flat $\g$-connections on $M$, then all characteristic
classes in $H^*_{\text{Lie}}(\g)$ are rigid.
\end{theorem}
\begin{proof}
The theorem follows from lemma \eqref{2.15} and theorem \eqref{3.6}.
\end{proof}

By checking dimensions in theorem \eqref{3.5}, exact sequence
\eqref{3.3}, and diagram \eqref{2.14}, we have:
\begin{corollary} \label{3.8}
If $HL^n(\g) = 0$ for $1 \leq n \leq p$, then all characteristic
classes in $H^n_{\text{Lie}}(\g)$ are rigid for $1 \leq n \leq p$.  
\end{corollary}

We close this section with two observations, one concerning a theorem
of P. Ntolo on the vanishing of $HL^* (\g)$ for $\g$ semi-simple, the
other concerning the highly nontrivial nature of $HL^* (W_1)$.
\begin{theorem} \label{3.9}
\cite{Ntolo}  If $\g$ is a semi-simple Lie algebra (over $\br$), then
$$  HL^n(\g) = 0 \ \ \ {\text{for}} \ \, n \geq 1.  $$
\end{theorem}
By contrast, the Leibniz cohomology of formal vector fields,
$HL^*(W_1)$, contains many non-zero classes which do not appear in
$H^*_{\text{Lie}}(W_1)$ \cite{Lodder}.  In the next section we compute
the image of a characteristic map
$$  HL^4(W_1) \to H^4_{\text{dR}}(M),  $$
where $M$ supports a family of codimension one foliations.

\section{Foliations}

Letting $W_1$ denote the Lie algebra of formal vector fields defined
in section two, recall that \cite[p.\ 101]{Fuks}
$$  H^{q-1}_{\text{Lie}}(W_1; \, W'_1) \simeq \br  $$
for $q=3$ and $q=4$, and zero otherwise.  
The generator of the class for $q=3$ is the
universal Godbillon-Vey invariant, called $\al$ in this paper, and we
denote the generator of the class for $q=4$ by $\zeta$.  From
\cite{Lodder}, the map
$$  H^{*-1}_{\text{Lie}}(W_1; \, W'_1) \to HL^*(W_1)  $$
given in exact sequence \eqref{3.3} is injective.  As dual Leibniz
algebras \cite{Bourbaki}, we have \cite{Lodder}
\begin{equation}  \label{4.1}
  HL^*(W_1 ) \simeq \Lambda(\al) \ot T(\zeta), 
\end{equation}
where $\Lambda (\al)$ is the exterior algebra on $\al$, and $T(\zeta)$
denotes the tensor algebra on $\zeta$.  

Let $M$ be a $C^{\infty}$ manifold with a one-parameter family ${\cal
F}_t$ of codimension one foliations having trivial normal bundles.
Let $\om_i(t)$ be the corresponding one-forms given in equation
\eqref{2.4} considered as differentiable functions of $t$.  Recall the
definitions of $\varPhi^*$ and $\phi_{\Theta}$ given in equations \eqref{2.12}
and \eqref{2.9} respectively.
In this section we prove the following:
\begin{theorem}  \label{4.2}
Let $M$ and ${\cal F}_t$ be given as above.  Then the composition 
$$  HL^4(W_1) \simeq H^3_{\text{Lie}}(W_1; \, W'_1) 
    \overset{\varPhi^*}{\rightarrow} H^4_{\text{Lie}}(\tilde{W}_1) 
    \overset{\phi_{\Theta}}{\rightarrow} H^4_{\text{dR}}(M)   $$
sends $\zeta$ to the de Rham cohomology class of 
$$ c(\zeta):= \om'_1(0) \wedge \om_0(0) \wedge \om_1(0) \wedge \om_2(0).$$
Moreover, the cohomology class of $c(\zeta)$ does not depend on the
choice\footnote{Since $c(\zeta)$ may also be written as 
$- \om'_1(0) \wedge \om_1(0) \wedge d\om_1(0)$, where $d$ denotes the
de Rham coboundary, it is not necessary to show that the class of $c(\zeta)$ is
independent of the choice of $\om_2(t)$.}
of $\om_0(t)$ or $\om_1(t)$.
\end{theorem}

\begin{proof}
We first compute $c(\zeta)$ on the level of cochains.  Consider the
vector space basis $\{ \be_i \}_{i \geq 0}$ of
${\text{Hom}}^c_{\br}(W_1, \, \br)$ given by
$$  \be_i \Big( \frac{x^j}{j!} \, \frac{d}{dx}\Big) = \delta_{ij}. $$
From \cite{Lodder}, the class of $\zeta$ in $HL^4 (W_1)$ is
represented by the cochain
$$  \be_1 \ot (\be_0 \wedge \be_1 \wedge \be_2).  $$
The cochain map
$$  i: \Omega^3(W_1; \, W'_1) \to C^4(W_1)  $$
inducing the isomorphism
$$  H^3_{\text{Lie}}(W_1; \, W'_1) \overset{\simeq}{\rightarrow}
	HL^4(W_1)  $$
satisfies 
$$  i\big( ( \be_0 \wedge \be_1 \wedge \be_2) \ot \be_1 \big) = 
	- \be_1 \ot ( \be_0 \wedge \be_1 \wedge \be_2).  $$
Also, it is know that
$$  ( \be_0 \wedge \be_1 \wedge \be_2) \ot \be_1  $$
generates $H^3_{\text{Lie}}(W_1; \, W'_1)$ (as an $\br$ vector space).

Let $v_1$, $v_2$, $v_3$, $v_4 \ \in T_p(M)$.  From the definition of
$\varPhi$ and $\phi_{\Theta}$, the image of $\zeta$ in $\Omega^4(M)$
is the 4-form which sends $v_1 \ot v_2 \ot v_3  \ot v_4$ to
\begin{align*}
- ( \be_1 \ot \be_0 \wedge \be_1 \wedge \be_2) \Big( & -A_1' \ot A_2
	\ot A_3 \ot A_4 + A'_2 \ot A_1 \ot A_3 \ot A_4  \\
  & -A'_3 \ot A_1 \ot A_2 \ot A_4 + A'_4 \ot A_1 \ot A_2 \ot A_3\Big),
\end{align*}
where
\begin{align*}
&  A_i = \sum_{n \geq 0} \, \om_n(0)(v_i) \, \frac{x^n}{n!} \,
\frac{d}{dx} \, , \ \ \ i = 1, \ 2, \ 3, \ 4,  \\
&  A'_i = \sum_{n \geq 0} \, \om'_n(0)(v_i) \, \frac{x^n}{n!} \,
\frac{d}{dx} \, , \ \ \ i = 1, \ 2, \ 3, \ 4.
\end{align*}
By the definition of the $\be_i$'s, the image of $\zeta$ is thus
$$ \big( \om'_1(0) \wedge \om_0(0) \wedge \om_1(0) \wedge \om_2(0)\big)
  (v_1 \ot v_2 \ot v_3 \ot v_4).  $$

To show that the de Rham cohomology class of $c(\zeta)$ does not
depend on the choice of $\om_0(t)$, consider the one-form
$$  u_0(t) = f \cdot \om_0(t),  $$
where $f:M \to \br$ is a $C^{\infty}$ function with $f(p) \neq 0$ for
all $p \in M$.  Letting $d$ denote the de Rham coboundary, we have
from equation \eqref{2.4}
$$  d \om_1(t) = \om_0(t) \wedge \om_2(t).  $$
Then
$$ \om'_0(0) \wedge \om_0(0) \wedge \om_1(0) \wedge \om_2(0) =
   - \om'_0(0) \wedge \om_1(0) \wedge d \om_1(0).  $$
Also,
\begin{align*}
& d u_0(t) = u_0(t) \wedge \Big( \frac{-df}{f} + \om_1(t) \Big) \\
& u_1(t) = \frac{-df}{f} + \om_1(t)  \\
& u'_1(t) \wedge u_1(t) \wedge du_1(t) = \om'_1(t) \wedge \om_1(t)
\wedge d \om_1(t) \\
& \ \ \ \ \ \ \ \ \ \ \ \ \ \ \ \ \ \ \ \ \ \ \ \ \ \ \ \ \ \ \ 
 + \frac{df}{f} \wedge \om'_1(t) \wedge d \om_1(t).
\end{align*}
It follows that 
\begin{align*}
  u'_1(0) \wedge u_1(0) \wedge du_1(0) = \; & \om'_1(0) \wedge \om_1(0)
	\wedge d \om_1(0) \\
& + d\big( \log ( |f| ) \, \om'_1(0) \wedge d \om_1(0) \big).
\end{align*}
Compare with Ghys \cite{Ghys}.  Of course,
$$ d(\om'_1(0)) = \om'_0(0) \wedge \om_2(0) + \om_0(0) \wedge
	\om'_2(0).  $$

To show that the cohomology class of $c(\zeta)$ does not depend on the
choice of $\om_1(t)$, consider the one-forms
$$  u(t) = \om_1(t) + f \cdot \om_0(t),  $$
where $ g: M \to \br$ is a $C^{\infty}$ function (which may have
zeroes on $M$).  Then
\begin{align*}
 u'(0) \wedge u(0) \wedge du(0)& = \om'_1(0) \wedge \om_1(0) \wedge d\om_1(0)\\
 & + g\cdot \om'_0(0) \wedge \om_1(0) \wedge d\om_1(0)   \\
 & + \om'_1(0) \wedge \om_1(0) \wedge dg \wedge \om_0(0) \\
 & + g \cdot \om'_0(0) \wedge \om_1(0) \wedge dg \wedge \om_0(0).
\end{align*}
It can be checked that
\begin{align*}
 & u'(0) \wedge u(0) \wedge du(0) = \om'_1(0) \wedge \om_1(0) \wedge
	d\om_1(0) + d(A), \\
 & A = g \cdot \om'_0(0) \wedge d\om_1(0) - dg \wedge \om'_0(0) \wedge
	\om_1(0) - \frac{1}{2} \, g^2 \cdot \om'_0(0) \wedge d \om_0(0).  
\end{align*}
\end{proof}

The paper is closed by noting that the current algebra $\tilde{\g}$ is
a Leibniz algebra in the sense of Loday \cite{Ens.Math.} with the
Leibniz bracket of $\varphi_1, \; \varphi_2 \in \tilde{\g}$ given by
$$  \langle \varphi_1(t), \, \varphi_2(t) \rangle =
    [\varphi_1(t), \, \varphi'_2 (0)]_{\text{Lie}},  $$
where $[ \ \, , \ ]_{\text{Lie}}$ is the usual Lie bracket on
$\tilde{\g}$, and $\varphi'_2(0)$ is the constant path at
$\varphi'_2(0)$.  The Leibniz bracket is
not necessarily skew-symmetric, 
$$  \langle \varphi_1(t), \, \varphi_2(t) \rangle \neq 
   - \langle \varphi_2(t), \, \varphi_1(t) \rangle ,  $$
but satisfies the following version of the Jacobi identity
$$ \langle \varphi_1(t), \, \langle \varphi_2(t), \, \varphi_3(t) \rangle \rangle 
  = \langle \langle \varphi_1(t), \, \varphi_2(t) \rangle, \,
  \varphi_3(t) \rangle - \langle \langle \varphi_1(t) , \,
  \varphi_3(t) \rangle, \, \varphi_2(t) \rangle,  $$
which is the defining relation for a Leibniz algebra.  Also see 
\cite{Balavoine} and \cite{Ens.Math.}.  

\bigskip
\noindent
ACKNOWLEDGEMENTS

\smallskip
\noindent
The author would like to express his gratitude to the Institut des
Hautes \'Etudes Scientifiques for their support while this paper was being
written.

\bigskip
\bigskip
\noindent
Math Sciences, Department 3MB

\noindent
New Mexico State University

\noindent
Las Cruces, NM  88003

\medskip
\noindent
e-mail:  jlodder@nmsu.edu


\begin{thebibliography}{99}

\bibitem{Balavoine} Balavoine, D., ``El\'ements de carr\'e nul dans
les alg\`ebres de Lie gradu\`ees,'' {\em Comptes Rendus Acad. Sci.},
S\'erie I, 319, (1994), 783--788. 

\bibitem{Fuks} Fuks, D.B., {\em Cohomology of Infinite-Dimensional Lie
Algebras}, (A.B. Sosinskii translator), Consultants Bureau, New York, 
1986.

\bibitem{Ghys} Ghys, E., ``L'invariant de Godbillon-Vey,'' {\em
S\'eminaire Bourbaki}, 706, (1988--89).  

\bibitem{Godbillon} Godbillon, C., Vey, J., ``Un invariant des
feuilletages de codimension 1,'' {\em Comptes Rendus Acad. Sci.},
S\'erie A, 273, (1971), 92--95.

\bibitem{Ens.Math.} Loday, J.-L., ``Une version non commutative des
alg\`ebres de Lie:  les alg\`ebres de Leibniz,'' {\em L'Enseignement
Math.}, 39, (1993), 269--293.

\bibitem{Bourbaki} Loday, J.-L., ``La Renaissance des Op\'erades,''
{\em S\'eminaire Bourbaki}, 792 (1994--95).  

\bibitem{overview}Loday, J.-L., ``Overview on Leibniz Algebras,
Dialgebras and Their Homology,'' {\em Fields Institute
Communications}, 17, (1997), 91--102.

\bibitem{LP} Loday, J.-L., Pirashvili, T., ``Universal Enveloping
Algebras of Leibniz Algebras and (Co)-homology,'' {\em Math. Annalen},
296 (1993), 139--158.

\bibitem{Lodder} Lodder, J.M., ``Leibniz Cohomology for Differentiable
Manifolds,'' {\em Annales Inst. Fourier, Grenoble}, 48, 1 (1998),
73--95.  

\bibitem{Ntolo} Ntolo, P., ``Homologie de Leibniz d'alg\'ebres de Lie
semi-simples,'' {\em Comptes Rendus Acad. Sci.}, S\'erie I, 318,
(1994), 707--710.  

\bibitem{Pirashvili} Pirashvili, T., ``On Leibniz Homology,'' {\em
Annales Inst. Fourier, Grenoble}, 44, 2, (1994), 401--411. 

\bibitem{Thurston}Thurston, W., ``Noncobordant Foliations of $S^3$,''
{\em Bulletin of the American Math. Soc.}, Vol. 78, No. 4, (1972),
511--514. 


\end{thebibliography}
\end{document}